\documentclass[a4paper,12pt]{article}
\usepackage{comment}
\usepackage{cite}
\usepackage{amsmath}
\usepackage{amssymb}
\usepackage{amsfonts}
\usepackage[T1]{fontenc}
\usepackage[utf8]{inputenc}
\usepackage{graphicx}
\usepackage{fancyhdr}
\usepackage{float}
\usepackage{xcolor}
\usepackage{authblk}
\usepackage{mathrsfs}
\usepackage{empheq}
\usepackage{url}

\usepackage{geometry}
\begin{document}

\title{Relations between $\pi$ and the golden ratio $\phi$ in the form of Bailey-Borwein-Plouffe-type formulas}

\author[$\dagger$]{Jean-Christophe {\sc Pain}\\
\small
CEA, DAM, DIF, F-91297 Arpajon, France\\
Université Paris-Saclay, CEA, Laboratoire Matière en Conditions Extr\^emes,\\ 
91680 Bruyères-le-Ch\^atel, France
}

\maketitle

\begin{abstract}
We provide a family of expressions of $\pi$ in terms of the golden ratio $\phi$ in the same spirit of the formula obtained by Bailey, Borwein and Plouffe for $\pi$. Connection with cyclotomic polynomials is outlined.
\end{abstract}

\section{Introduction}

In 1997, Bailey, Borwein and Plouffe published the following formula \cite{Bailey1997}:
\begin{equation}
\pi=\sum_{k=0}^{\infty}\frac{1}{16^k}\left(\frac{4}{8k+1}-\frac{2}{8k+4}-\frac{1}{8k+5}-\frac{1}{8k+6}\right),
\end{equation}
or
\begin{equation}
\pi=\sum_{k=0}^{\infty}\frac{1}{16^k}\frac{\left(47+151 k+120 k^2\right)}{\left(15+194 k+712 k^2+1024 k^3+512 k^4\right)},
\end{equation}
obtained using the parallel integer relation detection algorithm PSQL \cite{Bailey2001,Ferguson1999}. The strength of the latter formula is that it enables one to get the $k^{th}$ base-16 digit of $\pi$ without computing any other (prior) digit. Many other formulas of this kind, often referred to as ``BBP'' formulas, were obtained for the first powers of $\pi$ and other constants (Ap\'ery constant $\zeta(3)$, Catalan constant, $\ln 2$, \emph{etc.}) \cite{Bailey,Huvent2001,Adegoke2012,Adamchik1996,Adamchik1997,Chan2006,Chan2006b}. Finding exact mathematical relations between $\pi$ and the golden ratio
\begin{equation}
\phi=\frac{1+\sqrt{5}}{2}.
\end{equation}
is challenging \cite{Borwein1987,Valdebenito,Valdebenito2,Anderson2008,Adegoke2}. For instance, Baez obtained the following ``Vi\`ete-type'' formula \cite{Baez}:

\begin{equation}
\pi=\frac{5}{\phi}.\frac{2}{\sqrt{2+\sqrt{2+\phi}}}.\frac{2}{\sqrt{2+\sqrt{2+\sqrt{2+\phi}}}}.\frac{2}{\sqrt{2+\sqrt{2+\sqrt{2+\sqrt{2+\phi}}}}}\cdots
\end{equation}
%
and Chan found Machin-type formulas \cite{Chan2008,Luca2009}:
\begin{eqnarray}
\pi&=&4\arctan\left(\frac{1}{\phi}\right)+4\arctan\left(\frac{1}{\phi^3}\right)=8\arctan\left(\frac{1}{\phi^2}\right)+4\arctan\left(\frac{1}{\phi^6}\right)\nonumber\\
&=&12\arctan\left(\frac{1}{\phi^3}\right)+4\arctan\left(\frac{1}{\phi^5}\right).
\end{eqnarray}
Other interesting formulas were also recently published, derived from new forms of Taylor expansions of inverse tangent function \cite{Wu2022}. 

In the present  paper, we propose a family of BBP expressions for $\pi$ in terms of the golden ratio. The number $\pi$ and the golden ratio $\phi$ are related by
\begin{equation}
e^{i\frac{\pi}{5}}=\frac{\phi+i\sqrt{3-\phi}}{2},  
\end{equation}
which gives
\begin{equation}
\cos\left(\frac{\pi}{5}\right)=\frac{\phi}{2}
\end{equation}
and
\begin{equation}
\sin\left(\frac{\pi}{5}\right)=\frac{\sqrt{3-\phi}}{2}
\end{equation}
yielding
\begin{equation}
\tan\left(\frac{\pi}{5}\right)=\frac{\sqrt{3-\phi}}{\phi}
\end{equation}
which can be put in the integral form
\begin{equation}
\pi=5\int_0^{\frac{\sqrt{3-\phi}}{\phi}}\frac{dx}{1+x^2}.
\end{equation}
At this stage, it would be possible to consider the historical Madhava-Gregory-Leibniz, formula
\begin{equation}
\arctan(x)=\int\frac{dx}{1+x^2}=\sum_{k=0}^{\infty}\frac{(-1)^k}{2k+1}x^{2k+1}
\end{equation}
to get
\begin{equation}
\pi=\arctan\left(\frac{\sqrt{3-\phi}}{\phi}\right)=\sum_{k=0}^{\infty}\frac{(-1)^k}{2k+1}\left(\frac{3-\phi}{\phi^2}\right)^{k+1/2}.
\end{equation}
However, the approach we follow here, similar to the BBP one, is different. Since 
\begin{equation}
x^4-1=\left(x^2-1\right)\left(x^2+1\right),
\end{equation}
it is possible to write
\begin{equation}
\pi=5\int_0^{\frac{\sqrt{3-\phi}}{\phi}}\frac{1-x^2}{1-x^4}dx=5\sum_{k=0}^{\infty}\int_0^{\frac{\sqrt{3-\phi}}{\phi}}x^{4k}dx-5\sum_{k=0}^{\infty}\int_0^{\frac{\sqrt{3-\phi}}{\phi}}x^{4k+2}dx
\end{equation}
yielding
\begin{empheq}[box=\fbox]{align}
\pi=5\sum_{k=0}^{\infty}\left(\frac{\sqrt{3-\phi}}{\phi}\right)^{4k+1}\left[\frac{1}{4k+1}-\left(\frac{3-\phi}{\phi^2}\right)\frac{1}{4k+3}\right]%
\end{empheq}
In the same way, using
\begin{equation}
1-x^8=\left(1+x^2\right)\left(1-x^2+x^4-x^6\right),
\end{equation}
one gets
\begin{eqnarray}
\pi&=&5\int_0^{\frac{\sqrt{3-\phi}}{\phi}}\frac{1-x^2+x^4-x^6}{1-x^8}dx\nonumber\\
&=&5\sum_{k=0}^{\infty}\int_0^{\frac{\sqrt{3-\phi}}{\phi}}x^{8k}dx-5\sum_{k=0}^{\infty}\int_0^{\frac{\sqrt{3-\phi}}{\phi}}x^{8k+2}dx\nonumber\\
& &+5\sum_{k=0}^{\infty}\int_0^{\frac{\sqrt{3-\phi}}{\phi}}x^{8k+4}dx-5\sum_{k=0}^{\infty}\int_0^{\frac{\sqrt{3-\phi}}{\phi}}x^{8k+6}dx
\end{eqnarray}
yielding the final result:
\begin{empheq}[box=\fbox]{align}
\pi=&5\sum_{k=0}^{\infty}\left(\frac{\sqrt{3-\phi}}{\phi}\right)^{8k+1}\left[\frac{1}{8k+1}-\left(\frac{3-\phi}{\phi^2}\right)\frac{1}{8k+3}\right.\nonumber\\
&+\left.\left(\frac{3-\phi}{\phi^2}\right)^2\frac{1}{8k+5}-\left(\frac{3-\phi}{\phi^2}\right)^3\frac{1}{8k+7}\right].
\end{empheq}

\section{General case}

In the general case, for $p\in\mathbb{N}^*$:
\begin{equation}
x^{4p}-1=\left(x^2+1\right)\left(-1+x^2-x^4+\cdots + x^{2(2p-1)}\right),
\end{equation}
\emph{i.e.}
\begin{equation}
\frac{1}{x^2+1}=\frac{1}{\left(1-x^{4p}\right)}\sum_{k=0}^{2p-1}(-1)^{k}x^{2k},
\end{equation}
yielding the general formula, $\forall~p\in\mathbb{N}^*$:
\begin{empheq}[box=\fbox]{align}
\pi=5\sum_{k=0}^{\infty}\left(\frac{\sqrt{3-\phi}}{\phi}\right)^{4pk+1}~\sum_{i=0}^{2p-1}\left(\frac{\sqrt{3-\phi}}{\phi}\right)^{2i}\frac{(-1)^i}{4pk+2i+1},
\end{empheq}
which is the main result of the present work. Of course, since
\begin{equation}
\phi^2=\phi+1,
\end{equation}
one can write
\begin{empheq}[box=\fbox]{align}
\pi=5\sum_{k=0}^{\infty}\left(\frac{3-\phi}{\phi+1}\right)^{2pk+1/2}~\sum_{i=0}^{2p-1}\left(\frac{3-\phi}{\phi+1}\right)^{i}\frac{(-1)^i}{4pk+2i+1}.
\end{empheq}

\section{Conclusion}

We presented a family of expressions of $\pi$ in terms of the golden ratio $\phi$ in the same vein as the BBP formula derived for $\pi$. Although a number of relations involving simultaneously $\pi$ and $\phi$ exist in the literature, some of same being of the BBP type as well, the relations presented here were, to our knowledge, not published elsewhere. It is hoped that the procedure described in the present paper may stimulate the derivation of further relations. 

\section*{Appendix A: Connection with cyclotomic polynomials}

The cyclotomic polynomials can be defined as
\begin{equation}
\Phi_n(x)=\prod_{\substack{1\leq k\leq n\\ \gcd(k,n)=1}}\left(x-e^{\frac{2ik\pi}{n}}\right)
\end{equation}
and satisfy
\begin{equation}
\prod_{d|n}\Phi_d(x)=x^n-1.
\end{equation}
One has also
\begin{equation}
\Phi_n(x)=\prod_{d|n}\left(x^d-1\right)^{\mu\left(\frac{n}{d}\right)},
\end{equation}
where $\mu$ represents M\"obius function:
\begin{eqnarray}
\mu(n)&=&+1\;\;\;\;\mathrm{if}\;\; n\;\; \text{is a square-free positive integer with an even number of prime factors,}\nonumber\\
& &-1\;\;\;\;\mathrm{if}\;\; n\;\; \text{is a square-free positive integer with an even number of prime factors,}\nonumber\\
& &~~0\;\;\;\;\mathrm{if}\;\; n\;\; \text{has a squared prime factor.}
\end{eqnarray}
If $p$ is an odd prime number, one has
\begin{equation}
\frac{1}{1+x^2}=-\frac{\Phi_1(x)\Phi_2(x)\Phi_p(x)\Phi_{2p}(x)\Phi_{4p}(x)}{\left(1-x^{4p}\right)}
\end{equation}
where $\Phi_n$ are the cyclotomic polynomials \cite{Riesel1994}:
\begin{equation}
\Phi_1(x)=x-1,
\end{equation}
\begin{equation}
\Phi_2(x)=x+1,
\end{equation}
\begin{equation}
\Phi_p(x)=\frac{1-x^p}{1-x}=x^{p-1}+x^{p-2}+\cdots +x+1,
\end{equation}
\begin{equation}
\Phi_{2p}(x)=\frac{\left(1-x^{2p}\right)}{\left(1-x^p\right)}\frac{\left(1-x\right)}{\left(1-x^2\right)}=x^{p-1}-x^{p-2}+\cdots -x+1
\end{equation}
and
\begin{equation}
\Phi_{4p}(x)=\frac{\left(1-x^{4p}\right)}{\left(1-x^{2p}\right)}\frac{\left(1-x^2\right)}{\left(1-x^4\right)}=x^{2p-2}-x^{2p-4}+\cdots -x^2+1.
\end{equation}

\end{document}